\documentclass[12pt]{amsart}
\usepackage{amscd,amsmath,amssymb,amsfonts}
\usepackage[cmtip, all]{xy}
\theoremstyle{plain}

\theoremstyle{definition}

\numberwithin{thm}{section}
\numberwithin{equation}{section}

\newcommand{\ml}[2]{\begin{multline}\label{#1}#2 \end{multline}}
\newcommand{\ga}[2]{\begin{gather}\label{#1}#2 \end{gather}}

\newcommand{\surj}{\twoheadrightarrow}

\newcommand{\sF}{{\mathcal F}}

\newcommand{\sX}{{\mathcal X}}

\newcommand{\A}{{\mathbb A}}

\newcommand{\F}{{\mathbb F}}

\renewcommand{\P}{{\mathbb P}}
\newcommand{\Q}{{\mathbb Q}}

\newcommand{\X}{{\mathbb X}}

\newcommand{\Z}{{\mathbb Z}}

\begin{document}

\title[Higher congruences]{A remark on higher congruences for the number of rational points of varieties defined over a finite field}
\author{H\'el\`ene Esnault}
\address{
Universit\"at Duisburg-Essen,  Mathematik, 45117 Essen, Germany}
\email{esnault@uni-essen.de}

\date{July 16, 2005}
\begin{abstract}
We show that the $\ell$-adic cohomology of the mod $p$ reduction $Y$ of a regular model of a smooth proper variety defined over a local field, the cohomology of which is supported in codimension $\kappa$, can't be Tate up to level $(\kappa -1)$. As a consequence, the number of rational points of $Y$ can't fulfill the natural relation 
$|Y(\F_q)|\equiv \sum_{i\ge 0} q^i\cdot b_{2i}(\bar{Y}) $ modulo $q^\kappa$.
\\ \ \\

{\bf Une remarque sur les  congruences d'ordre sup\'erieur pour le nombre de points rationnels de vari\'et\'es d\'efinies sur un corps fini .}\\ \ \\
{\bf R\'esum\'e}: Nous montrons que la cohomologie $\ell$-adique de
la r\'eduction $Y$  modulo $p$ d'un mod\`ele r\'egulier d'une vari\'et\'e propre et lisse d\'efinie sur un corps local, dont la cohomologie est support\'ee en codimension $\kappa \ge 1$, ne peut \^etre de Tate jusqu'en niveau $(\kappa -1)$. En  cons\'equence, le nombre de points rationnels de $Y$ ne peut v\'erifier la formule naturelle $|Y(\F_q)|\equiv \sum_{i\ge 0} q^i\cdot b_{2i}(\bar{Y}) $ modulo $q^\kappa$. 

\end{abstract}
\maketitle
\begin{quote}

\end{quote}
{\bf Version fran\c{c}aise abr\'eg\'ee}. 
Dans  \cite{EsUneq}, Theorem 1.1, nous montrons que si $\sX$ est un mod\`ele r\'egulier d'une vari\'et\'e $X$ propre et lisse d\'efinie sur un corps local de corps r\'esiduel $\F_q$, alors si la cohomologie $\ell$-adique  $H^i(\bar{X})$ est support\'ee en codimension $\ge 1$ pour $i\ge 1$, le nombre de points rationnels de sa r\'eduction $Y$ modulo $p$ v\'erifie $|Y(\F_q)|\equiv 1 $ modulo $q$. 
En fait, pour \^etre plus pr\'ecis, sous cette hypoth\`ese, les valeurs propres du Frobenius g\'eom\'etrique agissant sur la cohomologie $\ell$-adique $H^i(\bar{Y})$ de $Y$ sont divisibles par $q$ en tant qu'entiers alg\'ebriques. Le but de cette note est de discuter une formulation en coniveau sup\'erieur. Une fa\c{c}on naturelle de g\'en\'eraliser la condition de coniveau $\ge 1$ pour $i\ge 1$ est de supposer que 
$H^i(\bar{X})/H^i(\bar{X})_{{\rm alg}}$ est support\'ee en codimension $\kappa$, o\`u $H^i(\bar{X})_{{\rm alg}}$ est nulle si $i$ est impair et sinon est la partie alg\'ebrique. Nous 
montrons cependant que cela n'implique pas que les valeurs propres du Frobenius g\'eom\'etrique sont divisibles par $q^\kappa$ en tant qu'entiers alg\'ebriques sur $H^i(\bar{Y})/H^i(\bar{Y})_{q^{\frac{i}{2}}}$, o\`u  $ H^i(\bar{Y})_{q^{\frac{i}{2}}}$ est nulle si $i$ est impair et sinon est la partie sur laquelle le Frobenius agit par multiplication par $q^{\frac{i}{2}}$. En particulier la formule naturelle 
$|Y(\F_q)|\equiv \sum_{i\ge 0} q^i\cdot b_{2i}(\bar{Y}) $ modulo $q^\kappa$ n'est pas valable en g\'en\'eral. Cette formulation est propos\'ee par N. Fakhruddin dans \cite{Fakh} qui la montre sous certaines hypoth\`eses pour une famille g\'eom\'etrique en \'egale caract\'eristique $p>0$. Nous montrons en quoi ces hypoth\`eses  sont tr\`es fortes. 

\section{Introduction}
In \cite{EsUneq}, Theorem 1.1, we show that if $\sX$ is a regular model of
a smooth proper variety $X$ defined over a local field
with finite residue field $\F_q$, then if $\ell$-adic cohomology 
$H^i(\bar{X})$ is supported in codimension $\ge 1$ for $i\ge 1$, the number of rational points of its mod $p$ reduction $Y$  fulfills $|Y(\F_q)|\equiv 1 $ modulo $q$. 
To be more precise, the assumption implies that the eigenvalues of the geometric Frobenius acting on $\ell$-adic cohomology $H^i(\bar{Y})$ of $Y$ are $q$-divisible algebraic integers. The proof relies on a version of Deligne's integrality theorem \cite{DeInt}, Corollaire 5.5.3 over local fields \cite{DE}, Corollary 0.4.
The goal of this note is to discuss a formulation in higher coniveau level. 
A natural generalization of the coniveau $\ge 1$ condition for $i\ge 1$ is to assume that $H^i(\bar{X})/H^i(\bar{X})_{{\rm alg}}$ is supported in codimension $\ge \kappa$, where $H^i(\bar{X})_{{\rm alg}}$ is equal to 0 if $i$ is odd, else is the algebraic part of cohomology. 
This means that there is a codim $\ge \kappa$ subscheme $Z\subset X$ so that $H^i(\bar{X})\xrightarrow{{\rm rest}=0} H^i(\bar{X}\setminus \bar{Z})/{\rm Im}(H^i(\bar{X})_{{\rm alg}}$. 
Said differently, $H^{i}(\bar{X})=H^{i}(\bar{X})_{{\rm alg}}$  for $i\le 2 \kappa$, and $H^i_{\bar{Z}}(\bar{X})\surj H^i(\bar{X})$ for $i\ge 2\kappa$.

However we show that this assumption does not imply that the eigenvalues of the geometric Frobenius acting on 
$H^i(\bar{Y})/H^i(\bar{Y})_{q^{\frac{i}{2}}}$ are divisible by $q^\kappa$-divisible algebraic integers,  where $ H^i(\bar{Y})_{q^{\frac{i}{2}}}$
is equal to 0 if $i$ is odd, else is the part of cohomology on which Frobenius acts by multiplication by $q^{\frac{i}{2}}$. In particular, the formula 
$|Y(\F_q)|\equiv \sum_{i\ge 0} q^i\cdot b_{2i}(\bar{Y}) $ modulo $q^\kappa$ does not hold in general. This formulation was proposed in  \cite{Fakh} by N. Fakhruddin, who shows  it under certain assumptions in a geometric family in equal characteristic $p>0$. We show how strong are those assumptions.

Our example consists  of a Godeaux surface in characteristic 0. We take a reduction mod $p$ which is a cone over a smooth curve $C$ of higher degree. After desingularization of the mod $p$ reduction, $H^1(\bar{C})(-1)$ enters $H^3(\bar{Y})$, and this destroys the possibility of the $|Y(\F_q)|\equiv 1 +q\cdot b_2(\bar{Y})$ mod $q^2$ congruence.

\noindent {\it Acknowledgements}. 
We thank Eckart Viehweg for discussions on the subject of this note and for his encouragement.

\section{The example}
Let us consider the Godeaux surface $X_0/\Q_p$ defined as the quotient
of the Fermat quintic $F\subset \P^3_{\Q_p}$ of homogeneous equation  $px_0^5+x_1^5+x_2^5+x_3^5$  by the group $\mu_5$ acting via $\xi\cdot (x_i)=(\xi^i\cdot x_i)$. Here $p$ is prime to 5, and $\xi$ generates the group of 5-th roots of unity. As well known \cite{BPVdV}, V, 15 and VII, 11, 
$H^0(X_0, \Omega^1_{X_0})=H^0(X_0, \Omega^2_{X_0})$ and by comparison of de Rham with \'etale cohomology,  one obtains $H^1(\bar{X}_0)=H^3(\bar{X}_0)=0, \ H^{2i}(\bar{X}_0)=H^{2i}_{{\rm alg}}(\bar{X}_0)$ for $i=0,1,2$. Let us assume we have a regular model $\sX\to {\rm Spec}(R)$
of $X_0$ over an extension $R\supset \Z_p$, with local field $K={\rm Frac}(R)$ and residue field $\F_q$. Thus the general fiber is $X=X_0\times_{\Q_p} K$, and we denote by $Y$ the mod $p$ reduction over $\F_q$. 

We use the computation in \cite{EsUneq}, sections 2 and 3. One has an exact sequence
\ga{2.1}{H^i_{\bar{Y}}(\sX^u) \to H^i(\bar{Y})\xrightarrow{{\rm sp}^u} H^i(X^u) \to H^{i+1}_{\bar{Y}}
(\sX^u)}
where $^u$ means the pull back via the extension $K\subset K^u$ to the maximal unramified extension, and $\bar{}$ means the pull back via the extension to the algebraic closure.  The sequence is equivariant with respect to the action of the geometric Frobenius ${\rm Frob}\in {\rm Gal}(\bar{\F}_q/\F_q) $ acting on 
$H^*(\bar{Y}), H^*_{\bar{Y}}(\sX^u), H^*(X^u)$. 
One also has the exact sequence
\ga{2.2}{0\to H^1(I, H^{i-1}(\bar{X}))\to H^i(X^u)\to H^i(\bar{X})^I\to 0}
where $I\subset {\rm Gal}(\bar{K}/K)$ is the inertia group, with quotient $
{\rm Gal}(\bar{K}/K)/I={\rm Gal}(\bar{\F}_q/\F_q)$. The sequence is equivariant with respect to the action of ${\rm Frob}$.
So using Gabber's purity theorem \cite{Fu}, Theorem 2.1.1. as in \cite{EsUneq}, section 2, one obtains
\ga{2.3}{H^1(\bar{Y})=0, }
and an equivariant exact sequence
\ga{2.4}{0\to H^0(\bar{Y}^0)(-1)\to H^2(\bar{Y}) \to H^2(\bar{X})^I,}
where $Y^0=Y\setminus $ singular locus. Thus in particular, Frob acts via multiplication by $q$ on $H^2(\bar{Y})$. 
So via Grothendieck-Lefschetz trace formula \cite{Gr} and the fact that $H^4(\bar{Y})=\oplus_{{\rm components}} \Q_\ell(-2)$, we conclude
\ga{2.5}{|Y(\F_q)| \equiv 1 +q\cdot b_2(\bar{Y}) - {\rm Tr}\ {\rm Frob}|H^3(\bar{Y}) \ {\rm mod} \ q^2.}
The question becomes  whether $H^3(\bar{Y})$ dies or not.

We now construct $\sX$ and show $H^3(\bar{Y})\neq 0$ for this $\sX$.  The mod $p$ reduction in $\P^3_{\F_p}$ of the model $\sF\subset \P^3_{\Z_p}$ of $F$  defined by the same equation  $px_0^5+x_1^5+x_2^5+x_3^5$ is the cone   over the Fermat curve $Q_{\F_p}\subset \P^2_{\F_p}$ of equation $x_1^5+x_2^5+x_3^5$. Then $\mu_5$ acquires one single fix point $(1:0:0:0) \in \P^3_{\F_p}$ which is the vertex of ${\rm cone}(Q_{\F_p})$. We base change $\Z_p\subset R$ via $\pi^5=p$ and denote by $k=\F_q$ the residue field and $K={\rm Frac}(R)\supset \Q_p$ the local field.
So $\sF\times_{\Z_p} R \subset \P^3_R$ is defined by the equation $\pi^5x_0^5+x_1^5+x_2^5+x_3^5$. The $\mu_5$ operation is still defined by $\xi\cdot x_i=\xi^i x_i$ and now the only fix point $x:=(1:0:0:0) \in \P^3_{\F_q}$ is at the same time the only point in which $\sF\times_{\Z_p} R$ is not regular. The affine equation of $\sF\times_{\Z_p} R$ in $(\A^3_R, x_0\neq 0)$ with coordinates $X_i=\frac{x_i}{x_0}$ on which  
$\mu_5$ acts via  $\xi\cdot X_i=\xi^iX_i$,  is $\pi^5+X_1^5+X_2^5+X_3^5$. 
We blow up the singularity $x$ to obtain $\sigma: \sF'\to \sF\times_{\Z_p} R$. Then $\sigma^{-1}(x)$ is isomorphic to the Fermat quintic $Z_2$ in $\P^3_{\F_q}$ of equation $X_0^5+X_1^5+X_2^5+X_3^5$ with action
$\xi\cdot X_i=\xi^i X_i$. Consequently, $\mu_5$ acts fix point free on $\sF'$ and 
the quotient $\sX$, which is defined over $R$, is a regular model of $X=X_0\times_{\Q_p} K:=(F/\mu_5)\times_{\Q_p} K$. Furthermore, $\sigma^{-1}(\sF\times_{\Z_p} \F_q)$ is the union of two components, one being the blow up $Z_1$ in the vertex of ${\rm cone}(Q_{\F_p}\times_{\F_p} \F_q)$, the other one being the Fermat quintic $Z_2$. Thus the mod $p$ fiber $Y$ of $\sX$ has two components
$S_i=Z_i/\mu_5$. They meet along $C=(Q_{\F_p}\times_{\F_p} \F_q)/\mu_5$.  As $p\ne 5$, the covering $Q_{\F_p}\times_{\F_p} \F_q \to C$ is \'etale, and 
${\rm genus}(C)=2$. 

The normalization sequence for $Y$ yields  a Frob  equivariant  exact sequence
\ga{2.6}{  H^3(\bar{Y})\to H^3(\bar{S}_1)\oplus H^3(\bar{S}_2)\to 0.}
On the other hand, one has 
\ml{2.7}{0\neq H^1(\bar{C})(-1)=H^1(\bar{Q}_{\F_p})^{\mu_5}(-1)=\\
H^3_c(\bar{Z}_1\setminus \bar{Q}_{\F_p})^{\mu_5}=H^3(\bar{Z}_1)^{\mu_5}=H^3(\bar{S}_1).} 
Thus
\ga{2.8}{H^3(\bar{Y})\surj H^1(\bar{C})(-1)\neq 0}
which shows $H^3(\bar{Y})\neq 0$.

\section{Discussion}
\subsection{Higher dimension} One can produce examples as above in all dimensions by taking the product $\sX\times_R \P^n$, which is still regular.  Then $H^i(X\times_K \P^n)/H^i(X\times_K \P^n)_{{\rm alg}}=0$ for all $i$, while $H^{3+2j}(Y\times_{\F_q} \P^n)\neq 0$ for all $j\ge 0$. 
\subsection{Motivic condition}  From  
\eqref{2.1},  using \eqref{2.2} and applying \cite{DE}, Corollary 0.4 to the eigenvalues of $H^i(X^u)$, we see immediately that the eigenvalues of Frob on $H^i(\bar{Y})$ fulfill
\ga{3.1}{{\rm sp}^u \ {\rm injective} \ + N^\kappa (H^*(\bar{X})/H^*(\bar{X})_{{\rm alg}})=
(H^*(\bar{X})/H^*(\bar{X})_{{\rm alg}})\\
\ \Longrightarrow {\rm eigenvalues \ Frob}|H^i(\bar{Y}) \notag \\ 
=\begin{cases}
0 & i< 2\kappa \ i \  {\rm odd}\\
q^{\frac{i}{2}} & i \le  2\kappa \  i \ {\rm even}\\
\in q^\kappa\cdot {\bar{\Z}} & i \ge 2\kappa. 
\end{cases}\notag
}
Here $N^\kappa$ is the coniveau filtration as explained in the Introduction.

In \cite{Fakh}, N. Fakhruddin analyzes the motivic conditions for a 
family $f: \X\to S$ defined over a finite field $k$, with $S, \X$ smooth, to have the property that a singular fiber $Y$ over a closed point $s$ with residue field $\F_q\supset k$ fulfills the property $|Y(\F_q)|=\sum_{i\ge 0} (-1)^i q^i\cdot b_{2i}(\bar{Y})$
modulo $q^\kappa$. More precisely, he studies the motivic conditions in a geometric family forcing the eigenvalue behavior described in \eqref{3.1}. He singles out three  conditions.
We explain them and analyze the consequences they have on the completion $\sX=\X\times_S R$ at $s$ of the family $f$. Here $R$ is the completion of the equal characteristic ring of functions at $s\in S$.
 Surely,  as in \cite{EsPoint}, the first one is  base change  for the Chow groups $CH_i(\bar{X}), i\le (\kappa -1)$. We know by Bloch's type argument that this implies the coniveau condition in level $\kappa$ on $H^*(\bar{X})/H^*(\bar{X})_{{\rm alg}}$, but we are extremely far of understanding that this is equivalent to it, as predicted by the general Bloch-Beilinson conjectures. The second one is that $R^if_*\Q_\ell$ are constant local systems. This is to say that the specialization map $H^i(\bar{Y})\to H^i(\bar{X})$ is an isomorphism, which in particular forces ${\rm sp}^u$ to be injective, but is stronger than this. So we see that those two conditions imply the weaker cohomological conditions in \eqref{3.1}
which already force the eigenvalue conclusion on $H^i(\bar{Y})$. The third condition says that the Chow groups $CH_i(\bar{Y}), i\le (\kappa -1)$, are hit by specialization. This should translate into the condition ${\rm sp}^u$ injective above, which is then a consequence
of the cohomological consequence of the condition forcing
$R^if_*\Q_\ell$ being a constant local system. 

At any rate, even if, as explained above,  the conditions developed in \cite{Fakh} are far from sharpness, they tacitly raise the question of a finer formulation, and are a motivation
for this note.

\subsection{Formula}
It is of course extremely rare that one can check motivic conditions. 
It is in the rule easier
to control cohomological conditions, and \eqref{3.1} gives conditions for a good behavior of rational points on $Y$. However, the condition ${\rm sp}^u$ injective is very nongeometric and likely very nonnatural as well. It would be better to understand a finer condition on the contribution of $H^i_{\bar{Y}}(\sX^u)$ in $H^i(\bar{Y})$ via the sequence \eqref{2.1}.

\bibliographystyle{plain}

\renewcommand\refname{References}

\end{document}